\documentclass{amsart}

\usepackage{amssymb} \usepackage{amsfonts} \usepackage{amsmath}
\usepackage{amsthm} \usepackage{epsfig} \usepackage{wrapfig}

\newcommand{\mettifig}[1]{\epsfig{file=#1}}

\newtheorem{lemma}{Lemma}
\newtheorem{thm}[lemma]{Theorem}
\newtheorem{prop}[lemma]{Proposition}
\newtheorem{cor}[lemma]{Corollary}

\theoremstyle{definition}
\newtheorem{defn}[lemma]{Definition}

\theoremstyle{remark}
\newtheorem{rem}[lemma]{Remark}

\newcommand{\matZ} {\ensuremath {\mathbb{Z}}}

\newcommand{\matRP} {\ensuremath {\mathbb{RP}}}

\author{Evgeny Fominykh}
\author{Bruno Martelli}

\address{Department of Mathematics \\%
Chelyabinsk State University \\%
ulitsa Brat'ev Kashirinykh 129 \\%
454021 Chelyabinsk, Russian Federation \newline
 \hspace*{10pt} Institute of Mathematics and Mechanics, Ural Branch of the Russian Academy of Sciences}

 \email{fominykh@csu.ru}
 \urladdr{http://www.topology.kb.csu.ru/$\sim$fominykh/}

\address{Dipartimento di Matematica \\%
Universit\`a di Pisa \\%
Via F.~Buonarroti 2 \\%
56127 Pisa, Italy%
}

 \email{martelli@dm.unipi.it}
 \urladdr{http://www.dm.unipi.it/$\sim$martelli/}

\title{$k$-normal surfaces}



\thanks{Both authors are partially supported by the INTAS
project ``CalcoMet-GT'' 03-51-3663. The first author
acknowledges the support of the Russian Foundation for Basic
Research (project 04-01-96014)}

\begin{document}

\begin{abstract}
Following Matveev, a \emph{$k$-normal surface} in a triangulated
$3$-manifold is a generalization of both normal and (octagonal) almost
normal surfaces. Using spines, complexity, and Turaev-Viro invariants of $3$-manifolds, we prove the following results:
\begin{itemize}
\item a minimal triangulation of a closed irreducible or a
bounded hyperbolic $3$-manifold contains no non-trivial $k$-normal sphere;
\item every triangulation of a closed manifold with at least $2$
  tetrahedra contains some non-trivial normal surface;
\item every manifold with boundary has only finitely many triangulations
  without non-trivial normal surfaces.
\end{itemize}
Here, triangulations of bounded manifolds are actually \emph{ideal} triangulations.
We also calculate the number of normal surfaces of
nonnegative Euler characteristics which are contained in the
conjecturally minimal triangulations of all lens
spaces $L_{p,q}$.
\end{abstract}

\maketitle

\section*{Introduction}
We study in this paper the existence and number of
(generalizations of) normal surfaces in various triangulated
$3$-manifolds. Theorem~\ref{main:thm} below
concerns the existence of (a generalization of) normal spheres in
minimal triangulations. Theorems~\ref{boundary:thm}
and~\ref{closed:thm} then show that, with finitely many exceptions,
all triangulations of a given manifold contain non-trivial normal
surfaces. Finally, Theorem~\ref{normlens2:thm} calculates the number
of normal surfaces of genus $0$ and $1$ in a family of triangulated
lens spaces.

\subsection*{Normal surfaces}
\begin{figure}
\mettifig{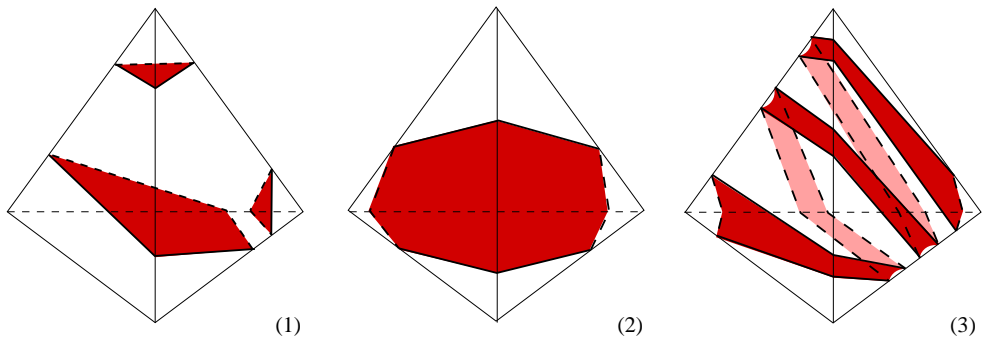, width = 13cm}
\caption{Some normal discs in a $k$-normal surface. In (3) we may wind the disc arbitrarily many times,
  getting discs with arbitrarily many intersections
  with the edges of the tetrahedron.
  }
\label{intersections:fig}
\end{figure}

A \emph{normal surface} in a triangulated $3$-manifold is a surface
intersecting each tetrahedron in triangles and squares, as in
Fig.~\ref{intersections:fig}-(1). Since every incompressible surface (in an
irreducible manifold) can be isotoped into normal position, normal
surfaces  can
be used to investigate the set of incompressible surfaces algorithmically.

An \emph{almost normal surface} is usually defined by admitting also
finitely many more types of intersections with some tetrahedra, such as
for instance the octagon of Fig.~\ref{intersections:fig}-(2). By a
beautiful result of Rubinstein and Thompson~\cite{Tho}, almost normal surfaces
can be used to find algorithmically some important types of compressible
surfaces, such as a $2$-sphere in $S^3$.

Matveev defined in~\cite{Mat_book} a \emph{$k$-normal surface} as a
surface that intersects each tetrahedron in polygons, having
consecutive vertices in distinct edges of the tetrahedron, and
incident to each such edge at
most $k$ times. With this definition, a normal surface is just $1$-normal, and an
almost normal surface with octagons is $2$-normal.
Polygons may have arbitrarily many edges as $k$ increases, as suggested by
Fig.~\ref{intersections:fig}-(3).

A surface which is $k$-normal for some $k$ was also called \emph{spinal} by Frohman and
Kania-Bartoszynska in~\cite{FroKa}. We prove here the
following.

\begin{thm} \label{main:thm}
Let $M$ be a closed irreducible or a bounded hyperbolic $3$-manifold,
different from $S^3, \matRP^3$, and $L_{3,1}$.
A minimal triangulation of $M$ does not contain any non-trivial
$k$-normal sphere, for all $k$.
\end{thm}

When the manifold is bounded here we actually
talk about \emph{ideal} triangulations,
see Section~\ref{triangulations:section} for definitions.
Theorem~\ref{main:thm} for normal spheres was known to Matveev in
1990, and its proof is essentially contained in~\cite[Theorem B]{Mat90}
as a byproduct of his \emph{complexity theory}. It was then extended
to the $2$-normal case by Jaco and Rubinstein \cite{JaRu}, see
\cite{Mat_book} for a short proof. We prove it in the general case
using spines. It turns out that $k$-normal surfaces
are sometimes better studied with spines rather than with triangulations:
the main reason for that is that by cutting a spine along a
$k$-normal surface we get another spine, while by cutting a triangulation we
get a polyhedral subdivision which is never a triangulation.

\subsection*{Existence of normal surfaces}
A connection between normal surfaces and Turaev-Viro invariants was
noted in \cite{FroKa}. We use here a nice version of the homology-free 5th
Turaev-Viro invariant, called the \emph{$t$-invariant}, due to
Matveev, Ovchinnikov, and Sokolov~\cite{t}, to prove the following results.

\begin{thm}\label{boundary:thm}
Every compact $3$-manifold has only finitely many triangulations without non-trivial normal surfaces.
\end{thm}
Actually, it turns out that such exceptional triangulations have a
fixed number of tetrahedra, depending only on the $t$-invariant of the
manifold. Despite this fact, there
are plenty of triangulations without non-trivial normal surfaces:
at least half of the $\sim 5000$ minimal triangulations of hyperbolic manifolds
with $\leqslant 4$ tetrahedra are so \cite{FriMaPe_pac, FriMaPe_exp}, and the number of
such triangulations with $n$ tetrahedra is bigger than $n^{cn}$, for
some $c>0$ and all $n>0$ \cite{FriMaPe_JDG}. (The
  proofs of \cite{FriMaPe_JDG} actually show the same results for triangulations
  without non-trivial $k$-normal surfaces.)

As the following result shows,
there is only one example with \emph{closed} manifolds, and it is the
minimal triangulation $T_{5,2}$ of the lens space $L_{5,2}$, having one vertex and one tetrahedron.
\begin{thm}\label{closed:thm}
Every triangulation of any closed manifold except $T_{5,2}$ contains a non-trivial normal surface.
\end{thm}

\begin{rem}
We actually prove in Section~\ref{constructions:section} a stronger
version of Theorems~\ref{boundary:thm} and~\ref{closed:thm}:
with the same exceptions, every triangulation contains
a non-trivial normal surface which intersects every edge in at most
$2$ points.
\end{rem}

The $t$-invariant is essential in the
proof of both theorems. Since their construction \cite{TuVi} in
1992, Turaev-Viro invariants have been used essentially only as a concrete
tool to quickly distinguish triangulated manifolds. This is the first ``theoretical'' application
of these invariants to $3$-dimensional topology known to the authors.

\subsection*{Triangulated lens spaces}
Using similar techniques we count the number of normal
surfaces of non-negative Euler characteristic in the conjecturally
minimal triangulations $T_{p,q}$ of all lens spaces. Let $S(p,q)$ be the sum of all partial quotients in the
expansion of $p/q$ as a regular continued fraction. The triangulation
$T_{p,q}$ contains $S(p,q)-3$ tetrahedra, see Section~\ref{lens:section}.

\begin{thm}
 \label{normlens2:thm}
 Let $(p,q)$ be coprime with $p\geqslant 4$.
 \begin{enumerate}
 \item $T_{p,q}$ contains no normal projective planes and
 no non-trivial normal spheres;
 \item the number $\tau (T_{p,q})$ of distinct (up to normal isotopy) normal tori in $T_{p,q}$ is
 $$\tau (T_{p,q})=\left\{ \begin{array}{rl}
  S(p,q)-3 & \mbox{if } q=1, \\  S(p,q)-4 & \mbox{in all other cases;} \end{array}
 \right.$$
\item the number $\kappa (T_{p,q})$ of distinct (up to normal isotopy) normal Klein bottles in $T_{p,q}$ is
 $$\kappa (T_{p,q})=\left\{ \begin{array}{rl} 1 & \mbox{if }
 p=4n \mbox{ and } q=2n\pm 1 \mbox{ for some positive integer } n, \\  0 & \mbox{in all other cases}. \end{array}
 \right.$$
 \end{enumerate}
\end{thm}

This paper is organized as follows. Precise definitions of all the
concepts introduce above are given in
Section~\ref{triangulations:section}. Theorem~\ref{main:thm} is proved
in Section~\ref{complexity:section}, while
Theorems~\ref{boundary:thm} and~\ref{closed:thm} are proved in
Section~\ref{constructions:section}. Theorem~\ref{normlens2:thm} is
proved in Section~\ref{lens:section}.

\subsection*{Acknowledgments} We warmthly thank Sergei Matveev for some helpful discussions.

\section{Surfaces and triangulations} \label{triangulations:section}
We fix here the terminology used in the Introduction and in the rest
of the paper.

Let $T$ be a simplicial face-pairing of some tetrahedra, and $|T|$ its support.
Here, $T$ is a \emph{triangulation} of a closed 3-manifold $M$ when $M = |T|$.
Tetrahedra are then not necessarily embedded in $M$. If $M$ has
boundary, we say that $T$ is a
\emph{triangulation} of $M$ if $M$ is $|T|$ minus open stars
of all its vertices (the latter is usually an \emph{ideal
  triangulation} in the literature: here we suppress the word
``ideal'' for simplicity).

Let now $T$ be a triangulation of some compact 3-manifold $M$. We say that a closed embedded
surface $\Sigma \subset M$ is \emph{$k$-normal} to $T$ if the following holds:
\begin{itemize}
\item $\Sigma$ is transverse to $T$ and disjoint from $\partial M$;
\item $\Sigma$ intersects each tetrahedron in discs, called
  \emph{normal discs};
\item $\Sigma$ intersects each triangular face in arcs connecting distinct edges;
\item each normal disc intersects each edge of the tetrahedron in at
  most $k$ points.
\end{itemize}
In a $1$-normal surface, normal discs are triangles and squares as in
Fig.~\ref{intersections:fig}-(1), while in a $2$-normal surface
they can also be octagons as in Fig.~\ref{intersections:fig}-(2). Normal discs can be
polygons with arbitrarily many edges as $k$ grows, as
suggested by Fig.~\ref{intersections:fig}-(3).

A normal surface is then just a $1$-normal surface. A \emph{trivial} normal surface is one made only of triangles.
Depending on whether the manifold is closed or not,
a connected trivial normal surface is either the link of a vertex or parallel to the boundary.

Two $k$-normal surfaces in a triangulation $T$, intersecting each edge of
$T$ in the same number of points are easily related via an isotopy,
called \emph{normal isotopy}.

A triangulation of a 3-manifold is \emph{minimal} if the manifold cannot be
triangulated with fewer tetrahedra.

We say that a compact $3$-manifold $M$ with boundary is \emph{hyperbolic} if it admits (after removing all tori
from $\partial M$) a complete hyperbolic metric of finite volume (possibly with cusps and geodesic
boundary). By Thurston's Hyperbolicity Theorem for Haken
  manifolds, this is equivalent to ask that
every component of $\partial M$ has $\chi\leqslant 0$, and $M$
does not contain essential surfaces with $\chi\geqslant 0$.

\section{Complexity} \label{complexity:section}
\subsection{Definition}
\begin{figure}
\mettifig{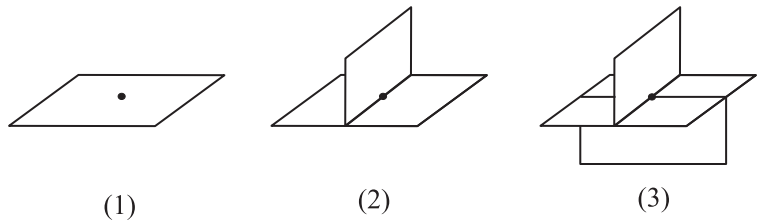, width = 10cm} \caption{Some neighborhoods of points
in an almost simple polyhedron.} \label{special:fig}
\end{figure}

The proof of Theorem~\ref{main:thm} relies heavily on complexity theory. We survey here
the definition and some of its properties. We start with the following.
\begin{defn}
A compact 2-dimensional polyhedron $P$ is \emph{almost simple} if
the link of every point in $P$ is contained in the graph
\mettifig{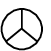, width = .4cm}.
\end{defn}
Alternatively, $P$ is almost simple if it is locally contained in
the polyhedron shown in Fig.~\ref{special:fig}-(3). A point, a
compact graph, a compact surface are therefore almost simple.
Three important possible kinds of
neighborhoods of points are shown in Fig.~\ref{special:fig}. A point
having the whole of \mettifig{mercedes_small.eps, width = .3cm} as a
link is called a \emph{vertex}, and its regular neighborhood is
shown in Fig.~\ref{special:fig}-(3). The set $V(P)$ of the vertices
of $P$ consists of isolated points, so it is finite. Note that
points, graphs, and surfaces do not contain vertices.

A compact polyhedron $P\subset M$ is a \emph{spine} of a compact manifold $M$ with boundary
if $M$ collapses onto $P$. When $M$ is closed, we say that $P\subset
M$ is a \emph{spine} of $M$ if $M\setminus P$ is an open ball.
\begin{defn}
The \emph{complexity} $c(M)$ of a compact connected 3-manifold $M$
is the minimal number of vertices of an almost simple spine of $M$.
\end{defn}
As an example, a point is a spine of $S^3$, and therefore $c(S^3)=0$.
The complexity of a
disconnected manifold is just defined as the sum of the complexities
of its components.

\subsection{Properties}
We list here some properties of $c$, proved
in~\cite{Mat90}, which will be used in the proof of
Theorem~\ref{main:thm}. The first one is the following.

\begin{thm} \label{addivity:thm}
Complexity is additive on connected sums. That is,
$$c(M\#M') = c(M) + c(M')$$
for all manifolds $M$, $M'$.
\end{thm}

\begin{figure}
\mettifig{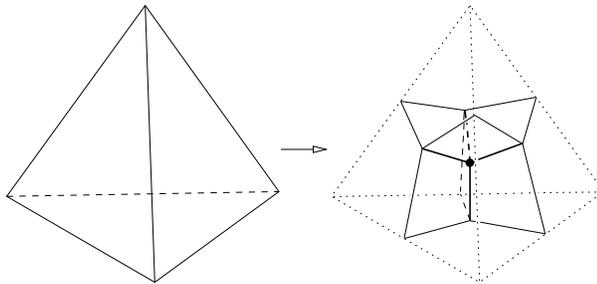, width = 8cm}
\caption{A spine of $M$ dual to a triangulation.}
\label{dualspine:fig}
\end{figure}

As suggested by Fig.~\ref{dualspine:fig},
a triangulation of a manifold $M$
with $n$ tetrahedra is dual to a polyhedron in $M$ with $n$
vertices. If $M$ is compact with boundary, the polyhedron is a spine of $M$. If it
is closed, $M\setminus P$ consists of some $h$ open balls, and by removing
open small discs from some appropriate $h-1$ faces of $P$ we get a
spine of $M$. In both cases, we get a spine with at most $n$ vertices.

This shows that $c(M)$ is less or equal than the number
of tetrahedra in a minimal triangulation of $M$.
One key result in complexity theory, due to Matveev~\cite{Mat90},
says that the converse is true on the most interesting $3$-manifolds.
\begin{thm} \label{triangulation:thm}
Let $M$ be a closed irreducible or a bounded hyperbolic $3$-manifold, different
from $S^3, \matRP^3$, and $L_{3,1}$. A minimal triangulation
of $M$ contains precisely $c(M)$ tetrahedra.
\end{thm}

\subsection{From $k$-normal surfaces to spines}
If $\Sigma\subset M$ is a closed embedded surface in a compact
$3$-manifold $M$, we denote by $M_\Sigma$ the (possibly
disconnected) compact manifold obtained from $M$ by removing an open
tubular neighborhood of $\Sigma$. We have the following general
result.

\begin{prop} \label{key:prop}
Let $T$ be a triangulation of a manifold $M$ with $n$ tetrahedra, and let
$\Sigma$ be a non-trivial $k$-normal surface. Then
$c(M_\Sigma)<n$.
\end{prop}

\begin{proof}
We extend the construction shown in
Fig.~\ref{dualspine:fig} to build a spine of $M_\Sigma$ with strictly
less than $n$ vertices.

Let $\Delta$ be an abstract tetrahedron of the triangulation. The surface
$\Sigma$ intersects $\Delta$ in discs, and each face of $\Delta$ into arcs connecting
distinct edges. For each such arc, draw another parallel dashed arc
closer to the nearest vertex, as in Fig.~\ref{construction:fig}-(1). Also draw a central $Y$-shaped
dashed graph as in the figure. Do this on each triangular face, so that the
dashed arcs and graphs match along their endpoints. The resulting
dashed object is a (possibly disconnected) graph $X$ in $\partial
\Delta$ with $4$
trivalent vertices.

\begin{figure}
\mettifig{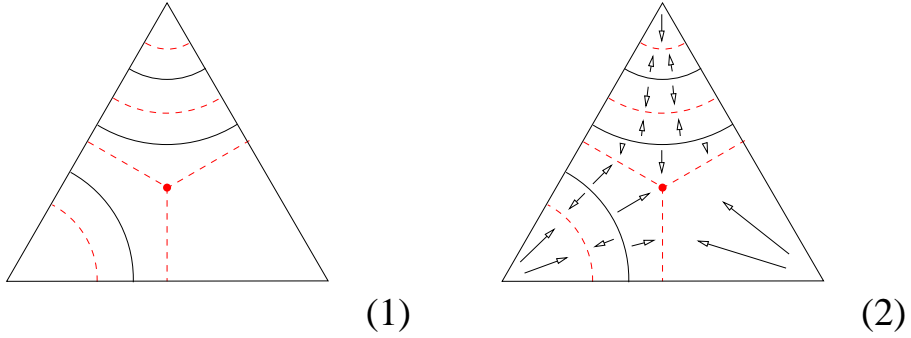, width = 12cm}
\caption{How to construct the graph $X$ (1), and a retraction onto it (2).}
\label{construction:fig}
\end{figure}

The complement of $\Sigma$ in $\Delta$ consists of
some topological balls. The boundary of each such abstract ball $B$
consists of some discs, arising from $\Sigma$, and one connected
component $\partial_v B$ coming from $\partial \Delta$, which contains a connected
component $X_0$ of the graph $X$. As shown in
Fig.~\ref{construction:fig}-(2), $X_0$ is a spine of
$\partial_v B\setminus\{$vertices of $\Delta\}$. We can extend the collapse to a collapse of $B$ onto
the subpolyhedron $c(X_0)\subset B$ obtained by coning
$X_0$ inside $B$.

We can force all dashed graphs and all collapses to match along the
edges and faces of the triangulation. As a result, the union $P$ of
all cones $c(X_0)$ is a spine of $M$. The polyhedron $P$ is made of
points as in Fig.~\ref{special:fig}-(1, 2), except possibly at the
vertices of the cones. Each $X_0$ is one of the following graphs:
\mettifig{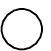, width = .3cm}, \mettifig{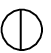, width =
.3cm}, \mettifig{mercedes_small.eps, width = .3cm}, and
\mettifig{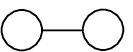, width = .8cm}. The cones on the first
three graphs yield the pieces shown in Fig.~\ref{special:fig}. The
cone on the fourth graph is not admitted in an almost simple spine;
however, it can be locally perturbed to an almost simple object as shown in
Fig.~\ref{perturbation:fig}.

\begin{figure}
\mettifig{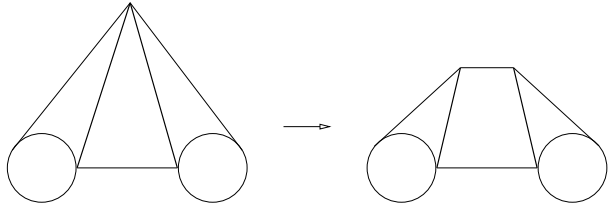, width = 8cm}
\caption{A perturbation of a cone into an almost simple polyhedron.}
\label{perturbation:fig}
\end{figure}

Therefore $P$ is (after a suitable perturbation) an almost simple
spine of $M$. Each \mettifig{mercedes_small.eps, width = .3 cm}
produces one vertex of $P$, while the other graphs yield no
vertices. There can be at most one such graph for each $\Delta$:
this gives $c(M_\Sigma)\leqslant n$. Equality holds if and only if
there is one  \mettifig{mercedes_small.eps, width = .3 cm} for each tetrahedron. This happens precisely
when $\Sigma$ intersects each tetrahedron into triangular faces,
\emph{i.e.}~when $\Sigma$ is trivial.

\end{proof}

\subsection{Proof of Theorem~\ref{main:thm}}
Let $M$ be closed irreducible or bounded
hyperbolic. A minimal triangulation $T$ contains precisely $c(M)$
tetrahedra, by Theorem \ref{triangulation:thm}.
If $\Sigma$ is a non-trivial $k$-normal sphere,
Proposition~\ref{key:prop} gives $c(M_\Sigma)<c(M)$. But $M$ is
irreducible, and hence $M_\Sigma = M\# B^3 \sqcup
B^3$; since $c(B^3)=0$ and $c$ is additive on connected sums, we get
$c(M_\Sigma) = c(M)$: a contradiction. Theorem~\ref{main:thm} is then proved.

\section{Constructions of normal surfaces} \label{constructions:section}

\subsection{Simple and special polyhedra}
To prove Theorems~\ref{boundary:thm} and~\ref{closed:thm} we need to introduce some
definitions.

A compact $2$-dimensional polyhedron is \emph{simple} if the link of
any point is homeomorphic to either \mettifig{circle.eps, width =
.3cm}, \mettifig{theta.eps, width = .3cm}, or
\mettifig{mercedes_small.eps, width = .3cm}. That is, its regular
neighborhood is one of those shown in Fig.~\ref{special:fig}. The
points of type (1), (2), and (3) form a stratification of $P$ into
vertices, \emph{$1$-components} and \emph{$2$-components}. The
$1$-stratum is the \emph{singular set} $S(P)$, consisting of all
points of type (2) and (3). If $1$-components are indeed open
segments and $2$-components are open discs, they are called
\emph{edges} and \emph{faces}, and the spine $P$ is \emph{special}.

As shown above in Fig.~\ref{dualspine:fig}, a
triangulation $T$ is dual to a special polyhedron $P$. When $M$ has
boundary, $P$ is a spine of $M$. If $M$ is closed, $P$ is a spine of
$M$ minus some open balls.

\subsection{Subpolyhedra and normal surfaces}
Let $P$ be a special spine. Assign a colour -- say, blue -- to some
faces of $P$. Let $Q$ be the union of the (closed) blue faces. It is
easy to see that:
\begin{itemize}
\item $Q$ is a simple polyhedron if and only if the number of blue
  faces incident to each edge of $P$ is either $0$, $2$, or $3$, with multiplicity;
\item $Q$ is a surface if and only if such a number is either $0$ or
  $2$.
\end{itemize}
Moreover, it is also easy to see that every simple subpolyhedron arises
in this way. This also shows that closed surfaces in $P$ are in natural
bijection with $H_2(P,\matZ_2)$.

\begin{figure}
\mettifig{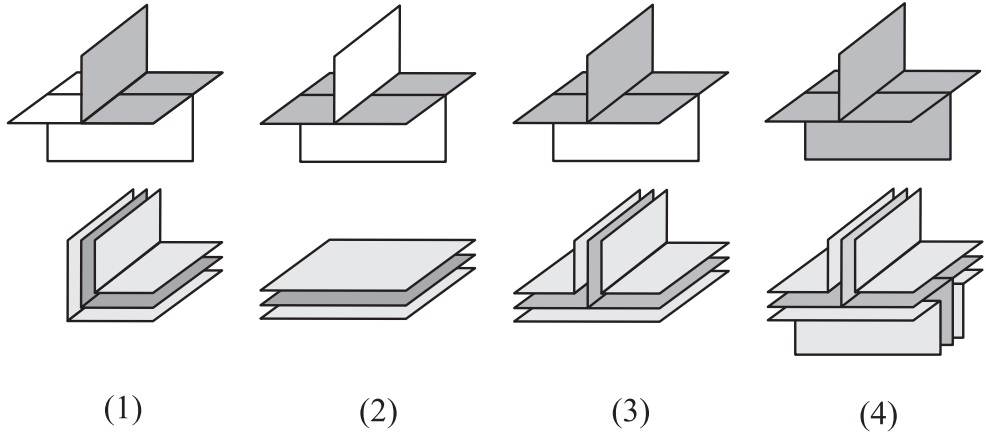, width = 13cm}
\caption{Four possible configurations of a simple subpolyhedron $Q$ near a
  vertex of $P$. The faces belonging to $Q$ are coloured. The boundary
  of a regular neighborhood consists respectively of (1) two triangles
  (2) two squares (3) two triangles and one square (4) four triangles.}
\label{simple:fig}
\end{figure}

Let $P$ be dual to a triangulation $T$.
The following facts are noted in~\cite{MatFom}:
\begin{itemize}
\item
if $Q$ is a surface, then it is normal in $T$;
\item
if $Q$ is a simple polyhedron, the boundary of a
small regular neighborhood of $Q$ is a normal surface in $T$.
\end{itemize}
A normal surface arising in this way is said to be respectively \emph{of
type I} and \emph{II}. In fact, if a vertex of $P$ belongs to $Q$, then
$4$ possible configurations may arise, as shown in
Fig.~\ref{simple:fig}: the same figure shows how each configuration
gives rise to normal triangles and squares.

\subsection{Universal subpolyhedron}
Let $T$ be a triangulation of a connected manifold $M$ with $k>0$ boundary components,
and let $P$ be the dual spine. Then $M\setminus P$ consists of
$k$ connected components. Assign a colour -- say, blue -- to those
faces of $P$ which belong to the closure of two distinct components
of $M\setminus P$. It is easy to see that the union $\Omega$ of the
(closed) blue faces is a simple subpolyhedron of $P$. The following
proposition shows that if $T$ does not contain non-trivial normal
surfaces of type I or II then $\Omega$ is in some sense universal.

\begin{prop} \label{only:surfaces:prop}
Suppose $T$ does not contain any non-trivial normal surface of type I or II.
Then $\Omega$ is a surface with $k-1$ connected components, each
parallel to some component of $\partial M$. Moreover, every proper
simple subpolyhedron of $P$ is a subsurface of $\Omega$.
\end{prop}
\begin{proof}
Suppose $P$ contains a proper simple subpolyhedron $Q$. As we said above, it
induces a normal surface of type I or II, which must hence be trivial.
Fig.~\ref{simple:fig} shows that configurations (2) and (3) give rise
to squares, and are thus ruled out. If one vertex of $P$ is as
in (4), all vertices of $P$ are so, because $P$ is connected: in this case $Q=P$. We are
therefore left only with vertices of type (1); in particular, $Q$ is a surface.

The union $S_1\cup S_2$ of two surfaces in $P$ is a simple
subpolyhedron, which is therefore either a surface or the whole of
$P$. The second case is easily ruled out since both $S_1$ and $S_2$
look as in Fig.~\ref{simple:fig}-(1) near the vertices. The
subpolyhedron $\Omega$ is in fact the union of $k$ surfaces, each consisting
of all faces separating one fixed boundary component to all
others. Therefore $\Omega$ is itself a surface (and not the whole of $P$).

Finally, every other proper simple
subpolyhedron $Q$ is a trivial (\emph{i.e.}~boundary-parallel)
surface, and hence contained in $\Omega$. It is also easy to
conclude that $\Omega$ has $k-1$ components.
\end{proof}

There are plenty of triangulations $T$ that satisfy the hypothesis of
Proposition~\ref{only:surfaces:prop}. As shown in~\cite{FriMaPe_JDG},
this holds for instance when $\partial M$ consists of a genus-$g$
surface and $h$ tori, and $T$ has $g+h$ tetrahedra. The number of such examples
with $n$ tetrahedra is bigger than $n^{cn}$ for some $c>0$ and all $n$, and the simplest one with at
least two boundary components has $3$ vertices \cite{FriMaPe_JDG}.

\subsection{Turaev-Viro invariants}
One of the simplest invariants of Turaev-Viro type is the
\emph{$t$-invariant} defined by Matveev, Ovchinnikov, and Sokolov in~\cite{t}. It equals the homology-free $5$-th Turaev-Viro invariant,
and is defined as follows. Let $\varepsilon$ be a fixed root of the polynomial
$x^2-x-1$. Let $P$ be a simple spine of a manifold $M$ with boundary. Then
$$t(P) = \sum_{Q\subset P} (-1)^{v_Q}\varepsilon^{\chi(Q)-v_Q}$$
where the sum is taken over all simple subpolyhedra $Q$ of $P$, and
$v_Q$ is the number of vertices of $Q$. It turns out that this
quantity does not depend on $P$ and gives therefore a well-defined invariant
$t(M)$ of $M$. If $M$ is a closed manifold, we define $t(M)$ as
$t(M\setminus B)$ where $B$ is a small open ball.

Like all (suitably normalized) Turaev-Viro invariants, $t$ is multiplicative on connected
sums (of closed manifolds) and $\partial$-connected sums (of manifolds
with boundary).

\subsection{Proof of Theorem~\ref{boundary:thm}}
Let $M$ be a manifold with $k\geqslant 1$ boundary components, and $T$ be a triangulation of
$M$ with $n$ tetrahedra. Let $P$ be its dual spine. We prove that,
with finitely many exceptions, $T$ contains a normal surface of type I
or II.

If not, by Proposition~\ref{only:surfaces:prop} there is a
boundary-parallel surface $\Omega\subset P$ containing all proper simple
subpolyhedra. Therefore
$$t(M) = t(P) = (-1)^n\varepsilon^{\chi(M)-n} + t(\Omega).$$
Suppose there are
infinitely many spines $P_1,\ldots,P_i,\ldots$ of $M$ without non-trivial normal
surfaces. Each $P_i$ has $n_i$ vertices and contains some universal $\Omega_i$
as above.

The quantity $(-1)^{n_i}\varepsilon^{\chi(M)-n_i} + t(\Omega_i)$
must not depend on $i$. Since $\Omega_i$ is boundary-parallel, it may vary only in a finite set of
topological types: hence $t(\Omega_i)$ varies in a finite set. On the
other hand, the number $n_i$ of vertices of $P_i$ tends to $\infty$,
since there are only finitely many special spines with a given number of vertices.
Therefore $(-1)^{n_i}\varepsilon^{\chi(M)-n_i}$ takes infinitely many
distinct values: a contradiction.

\subsection{Proof of Theorem~\ref{closed:thm}}
Let $T\neq T_{5,2}$ be a triangulation with $n$ tetrahedra of some closed
$3$-manifold $M$, and $P$ be the dual special spine. We prove that $T$
contains a non-trivial normal surface of type I or II.

If not, Proposition~\ref{only:surfaces:prop} implies that there is a
surface $\Omega\subset P$ containing all proper simple subpolyhedra of $P$.

Suppose first
that $T$ has $k\geqslant 2$ vertices. Therefore $P$ is a spine of $M$
minus $k$ small open balls. This manifold can be written as a
$\partial$-connected sum
$$(M\setminus B)\natural\underbrace{\big(S^2\times [0,1]\big)\natural\ldots\natural\big(S^2\times [0,1]\big)}_{k-1}$$
where $B$ is a small open ball in $M$.
By taking $S^2$ as a spine for $S^2\times [0,1]$ we see that $t\big(S^2\times
[0,1]\big) = \varepsilon^2+1$.
Since $t$ is multiplicative on $\partial$-connected sums, we have:
$$t(P) = t(M\setminus B) \cdot t\big(S^2\times [0,1]\big)^{k-1} = t(M)\cdot
(\varepsilon^2 + 1)^{k-1}.$$
On the other hand, $\Omega$ consists of $k-1$ disjoint
spheres, and therefore
$$t(P) = (-1)^n\varepsilon^{\chi(P)-n}+t(\Omega) = (-1)^n\varepsilon^{k-n} + (\varepsilon^2 + 1)^{k-1}.$$
The $t$-invariant is an element of the ring
$\matZ[\varepsilon]$. Since $k>1$, the equality
$$t(M)\cdot (\varepsilon^2 + 1)^{k-1} = (-1)^n\varepsilon^{k-n} +  (\varepsilon^2 + 1)^{k-1}$$
implies that $(\varepsilon^2 + 1)$ divides $\varepsilon^{k-n}$. It is
easy to see that $\varepsilon^2 + 1 =
\varepsilon + 2$ is not invertible in $\matZ[\varepsilon]$, while $\varepsilon^i$
is so for all $i$: a contradiction.

Suppose now that $T$ has one vertex. Let $f$ be a face of $P$.
By removing (the interior of) $f$ from $P$ we get an
almost simple polyhedron $P'$ such that $M\setminus P'$ is an open
solid torus. It collapses to a polyhedron $P''\subset P'$ which is the union of a simple polyhedron $Q$ and a
$1$-dimensional graph. (Every almost simple spine collapses to a
subpolyhedron of this kind, called \emph{nuclear} by Matveev in \cite{Mat_book}.)

Since $Q$ is a proper subpolyhedron of $P$, it is empty by our
assumption. Hence $P''$ is a
graph, and since $M\setminus P''$ is an open solid torus, the manifold $M$
is a lens space. It is then shown in
\cite[Theorem 8.1.28]{Mat_book} that $t(M)\in\{0,1,\varepsilon
+1,\varepsilon^2 + 1\}$. Since $P$ contains no proper subpolyhedron, we
have $t(P) = (-1)^n\varepsilon^{1-n}+1$, which belongs to $\{0,1,\varepsilon
+1,\varepsilon^2 + 1\}$ only if $n=1$ and $t(P)=0$.

Finally, a direct investigation shows that there is a unique
spine with $n=1$ vertices of a lens space with $t=0$. The lens space
is $L_{5,2}$ and the triangulation
is $T_{5,2}$, which indeed does
not contain any normal surface by Theorem~\ref{normlens1:thm} below, proved in \cite{MatFom}.

\begin{cor} Every special spine of a closed manifold except the
  minimal spine of $L_{5,2}$ contains a proper simple polyhedron.
\end{cor}

\section{Normal surfaces with nonnegative Euler
characteristics in lens spaces}\label{lens:section}
\subsection{The triangulations $T_{p,q}$.}
Every lens space $L_{p,q}$ has a candidate minimal triangulation $T_{p,q}$, which
is proved to be minimal only for finitely many values of $(p,q)$, and conjectured to
be so in all cases. We compute here the number of normal surfaces of
non-negative Euler characteristic in each such triangulation.

\begin{figure}
\mettifig{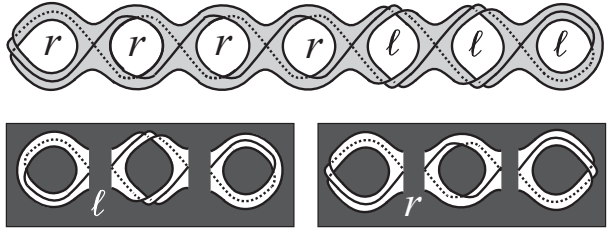, width = 8cm} \caption{How to construct a spine of a lens space.} \label{rlword:fig}
\end{figure}

We describe the triangulation $T_{p,q}$ by drawing regular
neighborhoods of the $1$-skeleton of the dual spine:
see~\cite{Mat_book} for details. Let $(p,q)$ be two positive coprime
number, with $p\geqslant 4$.
We define the operators $r,\ell
\colon Z\oplus Z \to Z\oplus Z$ via the rules $r(a,b)=(a,a+b)$ and
$\ell(a,b)=(a+b,b)$. The pair $(q,p-q)$ can be obtained from the
pair $(1,1)$ via a unique sequence $w$ of
operators $r, \ell$. We use the string $w$ to
construct the spine dual to $T_{p,q}$ as suggested in
Fig.~\ref{rlword:fig}-down. For example, if $(q,p-q)=(4,17)$, then $w=r r r
r \ell \ell \ell$. The spine of $L_{21,4}$
is shown in Fig.~\ref{rlword:fig}-up.

The following result was proved by Matveev and Fominykh
in~\cite{MatFom}.

\begin{thm} \label{FoMa}
 \label{normlens1:thm}
 Let $(p,q)$ be coprime with $p\geqslant 4$. Every connected normal
 surface in $T_{p,q}$
is either of type~I or of type~II.
\end{thm}

Let $S(p,q)$ be the sum of all partial quotients in the
expansion of $p/q$ as a regular continued fraction. The string $w$
contains $S(p,q)-2$ letters and thus $T_{p,q}$ contains $S(p,q)-3$ tetrahedra.

\subsection{Proof of Theorem~\ref{normlens2:thm}}.
Let $P$ be the special spine dual to $T_{p,q}$.
By Theorem~\ref{normlens1:thm}, every normal surface in
$T_{p,q}$ is determined by some simple subpolyhedron of $P$. Let us
analyze such polyhedra. One can easily see that every simple
subpolyhedron $G$ of $P$ contains the singular graph $S(P)$ of $P$.
This implies that $G$ can be obtained from $P$ by removing some
faces of $P$. Denote by $s$ the number of removed
faces. Since the Euler characteristic $\chi(P)$ of $P$ is
equal to $1$, we have $\chi(G) = 1-s$. Now we can conclude that
$T_{p,q}$ does not contain any normal projective plane and any
non-trivial normal sphere. Indeed, the boundary $F_G$ of a small
regular neighborhood of $G$, \emph{i.e.}~a normal surface of type II, is
orientable and its Euler characteristic is equal to $2-2s$.
Therefore, $F_G$ is a sphere if and only if $G=P$. But this sphere
is, obviously, trivial. Finally, since $P$ is not a surface, the
Euler characteristic of every normal surface of type I is less than
$1$.

As we know from~\cite{Mat_book}, $P$ has $S(p,q)-3$
vertices and $S(p,q)-2$ faces. Note that there are at most two faces
$f_1, f_2$ of $P$ that are incident more than once to some edge $e$.
It can be checked directly that $f_1 = f_2$ if and only if
the sequence $w$ determining $P$ consists of operators $r$ only,
\emph{i.e.}~if $q=1$. Since every simple subpolyhedron $G$ of $P$ contains
the singular graph of $P$, $G$ contains both $f_1$ and
$f_2$. Therefore, every normal surface of type I is
nonorientable.

 Let us prove the second assertion. As we already know, every normal
torus in $T_{p,q}$ is a surface of type II, \emph{i.e.}~it is the boundary
of a small regular neighborhood of some simple polyhedron $G$.
Moreover, $G$ is obtained from $P$ by removing a face
$f\neq f_1, f_2$ of $P$. If $p=5$ and $q=2$, then $P$
contains no proper simple subpolyhedra and hence $\tau(T_{5,2})=0$.
If $q=1$, then $f_1 = f_2$ and we have $S(p,q)-3$
possibilities to choose $f$. For all other pairs ($p$, $q$) we
have $f_1 \neq f_2$ and hence $\tau(T_{p,q})=S(p,q)-4$.

Now let $F$ be a normal Klein bottle. Then $F$ is a surface of type I.
Moreover, $F=P\setminus f$ for some face $f$ of
$P$. Since $F$ contains the singular graph of $P$, $f$ is the
unique face of $P$ incident once to each edge of $P$. Analyzing $P$ we
can directly see that such a face $f$ exists if and only
if the sequence $w$ determining $P$ is $w=r\underbrace{l\ldots
l}_{n-1}r$ for some positive integer $n$. Finally, we have $(2n-1,
2n+1) = r\underbrace{\ell\ldots \ell}_{n-1}r(1, 1)$ and hence $p=4n$,
$q=2n-1$.

Note that the conclusion (3) of Theorem~\ref{normlens2:thm}
agrees with the results of Rubinstein \cite{Rub}, since any Klein
bottle embedded in a lens space is incompressible and hence isotopic
to a normal surface.

\end{document}